\documentclass[a4paper,oneside,12pt]{article}
\usepackage{amsmath,amsfonts,amscd,amssymb}
\usepackage{longtable,geometry}
\usepackage[frenchb]{babel}
\usepackage[utf8]{inputenc}
\usepackage[T1]{fontenc}
\newtheorem{theorem}{Theorem}[section]
\newtheorem{corollary}[theorem]{Corollary}
\newtheorem{lemma}[theorem]{Lemma}
\newtheorem{proposition}[theorem]{Proposition}

\numberwithin{equation}{section}

\title{Law of the Iterated Logarithm for the random walk on the infinite percolation cluster}
\author{Hugo Duminil-Copin}
\date{September 2008}
\begin{document}
\maketitle
\small\textsc{Abstract:} We %%%
show that random walks 
on the infinite supercritical percolation clusters in $\mathbb{Z}^d$ satisfy the usual Law of the Iterated Logarithm. 
The proof combines Barlow's Gaussian heat
kernel estimates and the ergodicity of the random walk
on the environment viewed from the random walker as derived by Berger and Biskup. 
 \normalsize
\section{Introduction}

Asymptotic properties of random walks in $\mathbb {Z}^d$ are very well-understood. Their convergence to $d$-dimensional Brownian motions and their almost sure behavior (such as the law of the iterated logarithm) have been derived decades ago. A natural question to ask is what happens to random walks on graphs that are in some sense perturbations of $\mathbb{Z}^d$. One of the first examples to consider is to look at the random graph obtained by taking the infinite cluster $\mathcal {C}=\mathcal{C}(\omega)$ of a supercritical percolation process. One ``perturbs'' the original lattice by removing some edges independently. Various large-scale properties of this infinite graph have been studied with techniques such as coarse-graining. One of the most natural questions is to look at random walk on this cluster and to study its behavior.

One can for instance consider the {continuous-time simple random walk} (CTSRW) on $\mathcal{C}$. This is the process $X^{\omega}$ that waits an
exponential time of mean 1 at each vertex $x$ and jumps along one of the open edges $e$ adjacent to $x$, with each edge
chosen with equal probability. This process has been studied in a number of papers. Grimmett, Kesten, and Zhang
(\cite{GrimmettKestenZhang},1993) proved that $X^{\omega}$ is almost surely
recurrent if $d=2$ and transient if $d \geq 3$. Barlow (\cite{Barlow},2004)
proved Gaussian estimates for $X^{\omega}$. An invariance principle in every
dimension has been proved independently by Berger and Biskup in
(\cite{BergerBiskup},2004) and by Mathieu and Piatnitski
(\cite{MathieuPiatnitski},2004). Before that, Sidoravicius and Sznitman proved this result
for $d\geq 4$ (\cite{SidoraviciusSznitman}, 2004).
All these results show that a property that holds for random walk on $\mathbb{Z}^d$ still holds for random walk
on the infinite supercritical 
percolation cluster. 

It is natural to ask if this is still valid if one looks for instance at almost sure properties of the random walk 
(recall that almost sure properties often describe the behavior of the walk at exceptional times). 
Our goal in the present note is to show that it is indeed the case for the law of the iterated logarithm (LIL).

\begin{theorem} \label{law of the iterated logarithm} 
Consider $d \ge 2$ and suppose that $p>p_c$, where $p_c= p_c (d)$ is the critical bond percolation probability in $\mathbb {Z}^d$.
Then, there exists a positive and finite constant $c(p,d)$, such that for almost all realization of percolation with parameter $p$,  for all $x$ in the infinite cluster $\mathcal{C}$, the continuous-time random walk $X^{w}$ started from $x$  satisfies almost surely the following LIL:
$$\limsup_{t\rightarrow \infty} \frac{\left|X_t^{\omega}\right|}{\sqrt{t\log \log t}}=c(p,d).$$
\end{theorem}

Here and throughtout the paper, $|x| = |x|_1 = \sum_{j=1}^d |x_j|$ stands for the $L^1$ norm of  $x = (x_1, \ldots, x_d) \in {\mathbb Z}^d$.
Our proof can trivially be adapted to other norms (this would just change the value of the constant).
Note also that we are studying almost sure properties of the walk, so that the annealed and quenched statements are identical here (once we say that the constant $c(p,d)$ does not depend on the environment). 

The main ingredients of our proof are the Gaussian bounds
derived in  Barlow \cite{Barlow} and the ergodicity of Kipnis-Varadhan's \cite {KipnisVaradhan}  random walk on the environment as seen from the random walker derived by Berger and Biskup in  \cite{BergerBiskup}. These two results have in fact been instrumental in the (much more difficult) derivation of  the invariance principle for this random walk.

The paper is organized as follows. In Section 2, we will show that one can find positive and finite $c_1 (p,d)$ and $c_2 (p,d)$ such that almost surely
$$ c_1(p,d) \le \limsup_{t\rightarrow \infty} \frac{\left|X_t^{\omega}\right|}{\sqrt{t\log \log t}} \le c_2(p,d).$$
This will be based on the Gaussian estimates derived in \cite{Barlow}. 
%Note that one can not obtain ``uniform'' Gaussian bounds on
%$\mathbb{P}_{\omega}(X^{\omega}_t=y)$ for every $x,y\in \mathcal{C}$ and $t\geq 0$, because the local
%geometry of the infinite cluster prevents us from properly estimating the heat kernel for small $t$. However, 
%since we are studying an asymptotic property, it will still be possible to exploit these Gaussian bounds for large $t$. 
The upper bound is an easy application of the Borel-Cantelli Lemma whereas the proof of the lower bound will use the Markov property and the fact that  one can apply Gaussian bounds uniformly for $x$ in a ball of sufficiently large radius depending on $t$.
 
In Section 3, we derive a Zero-one Law for the limit of a discrete analog of the CTSRW. The main ingredient will be the ergodicity of a certain 
 shift $T$, related to Kipnis-Varadhan's {random walk on the
environment}. It has been proved by Noam Berger and Marek Biskup \cite{BergerBiskup} that this shift $T$ is
ergodic. Translating properties of the random walk in terms of this shift will allow us to
derive a Zero-one law for the limsup in the LIL for this discrete-time random walk. 

Finally, in Section 4, we conclude by checking that the time-scales of  the discrete-time random walk and of the continuous-time random walk are comparable.

\section{Weak LIL for the continuous random walk}

We will consider Bernoulli bond percolation of parameter $p$ on $\mathbb{Z}^d$ defined on a probability space
($\Omega$,$\mathcal{F}$,$\mathbb{P}_p$). It is well known (Grimmett \cite{grimmett}) that there exists $p_c\in(0,1)$ such
that when $p>p_c$ there is a unique infinite open cluster, that we denote by $\mathcal{C}$. For ${\mathbb P}_p$ almost every
environment $\omega \in \Omega$ and $x \in \mathcal {C}$, we define a CTSRW $X^{\omega}=(X_t^{\omega},t\geq0)$ started from $x$ under the probability measure $\mathbb{P}_{\omega}^x$. \textbf{In the whole paper, we fix $p>p_c$ and $d\geq2$}.

Because of translation-invariance of our problem (and because we are dealing with almost sure properties), 
we can restrict ourselves to the case $x=0$ and work with the probability measure 
 $\tilde{\mathbb{P}}_p=\mathbb{P}_p(.\left|0 \in \mathcal{C}\right.)$ and $X^{\omega}_0=0$.
We will use the notation $\Phi(t)=\sqrt{t\log \log t}$ for all $t>e$. 

\medbreak

We now recall Barlow's Gaussian estimates.
 The first one uses the 
{chemical distance} $d_{\omega}$ (or graph distance) on
$\mathcal{C}$. For every $x$ and $y$ in
$\mathcal{C}$, $d_{\omega}(x,y)$ is the length of the shortest
path between $x$ and $y$ that uses only edges in $\mathcal{C}$. For every integer $n$ and $x\in \mathcal{C}$, $\mathcal{B}_{\omega}(x,n)$ will denote the ball of radius $n$ and of center $x$ for distance $d_{\omega}$.

\begin{proposition} \label{maximum}(Barlow, \cite{Barlow}) There exist two constants $a_1=a_1(p,d)$ and $a_2=a_2(p,d)$ such that for every $\gamma>0$, there exists a finite random variable
$M_{\gamma}$ satisfying for almost every environment
$\omega$:
$$\text{for all } n\geq M_{\gamma}(\omega),\ \mathbb{P}^0_{\omega}(\max_{k\in \left[0,n\right]} d_{\omega} (0,X_k^{\omega}) > \gamma
\Phi(n)) \leq a_1 \exp\left(-a_2\frac{(\gamma\Phi(n))^2}{n}\right).$$
\end{proposition}
Recall that this statement holds for a general class of graphs (see Proposition 3.7 of Barlow \cite{Barlow}); percolation
estimates (see Theorem 2.18 and Lemma 2.19 of Barlow \cite{Barlow}) show that the percolation cluster belongs to this class. 
The other result that we will use is the Gaussian bound itself: 
\begin{theorem} (Barlow, \cite{Barlow}) \label{Gaussian bound} There exist finite constants $c_1$,..., $c_8$ and $\epsilon>0$ only depending on $p$ and $d$ that satisfy the following property. There exists a random variable $S_0$ with $\tilde{\mathbb{P}}_p(S_0\geq n)\leq c_7\exp(c_8 n^{\epsilon})$ and for almost every environment
$\omega$ such that $0,y\in \mathcal{C},t\geq 1$: 

(1) The transition density $p_t^{\omega}(0,y)$ of $X^{\omega}$ satisfies the Gaussian bound
$$c_1t^{-d/2}e^{-c_2{\left|y\right|^2}/{t}} \leq p_t^{\omega}(0,y) \leq c_3t^{-d/2}e^{-c_4{\left|y\right|^2}/{t}}\text{    for }t\geq S_0(\omega)\vee \left|y\right|.$$

(2) $c_5n^d \leq Vol(\mathcal{B}_{\omega}(0,n))\leq c_6 n^d\ \text{for}\ n\geq
S_0(\omega).$
\end{theorem}

Note that translation invariance makes possible for each $x \in \mathbb{Z}^d$, a random variable $S_x$
satisfying the analogous conditions (with the same constants $c_1$,...,$c_8$,$\epsilon$) where one just replaces the origin $0$ by $x$ (and therefore replaces $y$ by $x+y$).

Let remark that there is no uniform Gaussian bounds for every $x,y\in \mathcal{C}$ and every $t>0$ because (almost surely) every finite graph is actually embedded somewhere in the infinite cluster. We can now derive almost sure upper and lower bounds for our limsup.

\begin{proposition} \label{upper bound} \textbf{(Upper bound)} There exists a finite $c_+=c_+(p,d)$ such that for almost
every environment $\omega$,
$$\mathbb{P}^0_{\omega}\ a.s.\ \limsup_{t\rightarrow \infty} \frac{\left|X_t^{\omega}\right|}{\Phi(t)}\leq c_+.$$
\end{proposition}

\textbf{Proof:} Fix $\omega$ an environment containing 0. The proof goes along the same lines as in the Brownian case. Let $\gamma>0$, and define the following events:
$$A_n^{\omega}=\left\{\max_{k \in \left[0,2^n\right]}d_{\omega}(0,X_k^{\omega}) > \gamma \Phi(2^n) \right\}.$$
Proposition \ref{maximum} shows that for all $n$ large enough, 
$$\mathbb{P}^0_{\omega}(A_n^{\omega})\leq a_1\exp\left(-a_2\frac{(\gamma\Phi(2^n))^2}{2^n}\right)\leq 2a_1 n^{-a_2\gamma^2}.$$
Providing $\gamma$ large enough, the Borel-Cantelli Lemma claims that almost surely $A_n^{\omega}$ holds finitely often. Using the fact that $\left|.\right| \leq d_{\omega}(0,.)$, we get that for $n$ large enough, $\max_{k\in[0,2^n]}\left|X_k^{\omega}\right| < \gamma \Phi(2^{n})$. We conclude that for $n$ large enough, $\left|X_n^{\omega}\right| < 2\gamma \Phi(n)$.
%\qed
\begin{flushright}
$\square$
\end{flushright}

\begin{proposition} \label{lower bound} \textbf{(Lower bound)} There exists a positive $c_-=c_-(p,d)$ such that for almost
every environment $\omega$,
$$\mathbb{P}^0_{\omega}\ \text{a.s.},\ c_-\leq \limsup_{t\rightarrow \infty} \frac{\left|X_t^{\omega}\right|}{\Phi(t)}.$$
\end{proposition}

Let first present the outline of the proof. Consider $q>1$ and $\gamma>0$ (we will choose their values later). As in the Brownian case, set
$D_{n}^{\omega}=X_{q^n}^{\omega}-X_{q^{n-1}}^{\omega}$. We have $\left|X_{q^n}^{\omega}\right|\geq
\left|D_n^{\omega}\right|-|X_{q^{n-1}}^{\omega}|$. Using the upper bound, we obtain that almost surely, for $n$ large enough:
\begin{align}\left|X_{q^n}^{\omega}\right|&\geq \left|D_n^{\omega}\right|-2c_+\Phi(q^{n-1}).\end{align}Because $\Phi(q^{n-1})\leq q^{-1/2}\Phi(q^n)$, the second term can be chosen much smaller than $\Phi(q^n)$, providing $q$ large enough. Then, in order to prove the result, it is enough to bound  $D_n^{\omega}$ from below. Define the events $C_n^{\omega}=\left\{\left|D_n^{\omega}\right|>\gamma \Phi(q^n)\right\}.$ If these events hold for infinity many $n$ almost surely, then we are done. We define the $\sigma$-fields $\mathcal{F}_{n}^{\omega}=\sigma(X^{\omega}_k,k\leq q^{n})$. We will apply
the Borel-Cantelli Lemma generalized to dependent events (see Durrett \cite{Durrett}, chapter 4, paragraph 4.3). We therefore need to prove that
$$\mathbb{P}_{\omega}^0\text{ a.s. }\sum_{n\geq 1} \mathbb{E}^0_{\omega}[C_n^{\omega}
\left|\mathcal{F}^{\omega}_{n-1}\right.]=\infty.$$
Using the Markov property and Gaussian bounds, we will be able to 
find a lower bound for  $\mathbb{E}^0_{\omega}[C_n^{\omega}
\left|\mathcal{F}^{\omega}_{n-1}\right.].$ 
In order to apply these bounds, we need to control not only $S_0$ (from Theorem \ref{Gaussian bound}) but also $S_x$ for $x=X^{\omega}_{q^{n-1}}$. We first prove that it is indeed possible, using Gaussian estimates and the upper bound.

\begin{lemma} \label{environment}
Let $\gamma>0$, for almost every environment $\omega$ we have almost surely $S_{X^{\omega}_{n}}\leq \gamma\Phi(n)$ for $n$ large enough.
\end{lemma}

\textbf{Proof:} Let $\gamma>1$. Define for each integer $n$ the set 
$$B_n=\left\{\exists y\in B(0, 2c_+\Phi(n))\text{ s.t. } S_y\geq \gamma\Phi(n)\right\}.$$ 
where $S_y$ is the random variable of Theorem \ref{Gaussian bound}. The Theorem yields
$$\tilde{\mathbb{P}}_p(B_n)\leq Vol\left(B(0,2c_+\Phi(n))\right) d_1 \exp\left(-d_2(\gamma\Phi(n))^{\epsilon}\right).$$
The right-hand side of the inequality is summable, so that (by Borel-Cantelli) $B_n$ holds for a finite number of $n$  for almost every environment. But almost surely, $X_{n}^{\omega}$ is less than $2c_+\Phi(n)$ for $n$ large enough. Combining these two facts, we obtain the claim.
\begin{flushright}
$\ $
\end{flushright}
$ $\\
\textbf{Proof of Proposition \ref{lower bound}:} Let $q,\gamma>0$ and $\kappa>0$ such that $c_5\kappa^d>c_6+1$. Note that $\kappa$ does not depend on $\gamma$ and $q$. Set $t_n=q^n-q^{n-1}$. By the Markov property, we get for $n\geq1$,$$\mathbb{E}^0_{\omega}(C_n^{\omega} \left|\mathcal{F}^{\omega}_{n-1}\right.)=\mathbb{P}^0_{\omega_n}[\gamma \Phi(q^n)<X_{t_n}^{\omega_n}]\geq \mathbb{P}^0_{\omega_n}[\gamma \Phi(q^n)<X_{t_n}^{\omega_n}<\kappa\gamma \Phi(q^n)]=G_n(\omega_n)$$ 
where $\omega_n=\tau_{X_{q^{n-1}}^{\omega}}(\omega)$ ($\tau_x$ is the shift defined by $(\tau_x\omega)_y=\omega_{x+y}$) and:$$G_n(\omega)=\mathbb{P}^0_{\omega}\left[\gamma\Phi(q^n)<X^{\omega}_{t_n}<\kappa\gamma\Phi(q^n)\right].$$
The function $G_n$ is well-defined and measurable. If $\mathcal{A}_n(\omega)$ is the annulus $$\mathcal{A}_n(\omega)=\left\{z\in \mathcal{C},\text{ s.t. }\gamma\Phi(q^n)< \left|z\right| <\kappa \gamma \Phi(q^n) \right\},$$ we find by definition of the transition density $G_n(\omega)=\sum_{z\in \mathcal{A}_n(\omega)}p_{t_n}^{\omega}(0,z)$. We deduce:
\begin{align}\mathbb{E}^0_{\omega}(C_n^{\omega} \left|\mathcal{F}^{\omega}_{n-1}\right.)=\sum_{z\in \mathcal{A}_n(\omega_n)}p_{t_n}^{\omega_n}(0,z).\end{align}
Using Lemma \ref{environment}, we know that almost surely there exists $N$ large enough such that for every $n$ larger than $N$, $S_0(\omega_n)=S_{X^{\omega}_{q^{n-1}}}(\omega)\leq \gamma\Phi(q^{n-1})\leq t_n$. For $n\geq N$, one can use Gaussian estimates of Theorem \ref{Gaussian bound} for every $z\in \mathcal{A}_n(\omega_n)$, we get for such a $z$:
\begin{align}p_{t_n}^{\omega_n}(0,z)&\geq c_1t_n^{-d/2}\exp\left(-\frac{c_2(\kappa\gamma\Phi(q^n))^2}{t_n}\right)\geq c_1t_n^{-d/2}n^{-c_2(\kappa\gamma)^2}\end{align}
Using again the same Lemma, Theorem \ref{Gaussian bound} yields that the volume growth property holds for $S_n(\omega_n)$. Recalling the definition of $\kappa$, we find:
\begin{align}Vol(\mathcal{A}_n(\omega_n))&\geq (\gamma \Phi(q^n))^{d}\geq \gamma^d t_n^{d/2}\end{align}
Combining (2.3) and (2.4) in (2.2), we obtain that there exists a constant $c>0$ such that almost surely for $n$ large enough:
$$G_n(\omega_n)\geq cn^{-c_2(\kappa\gamma)^2}$$
Providing $\gamma$ small enough, we can use the generalized Borel-Cantelli Lemma (e.g. \cite{Durrett}). We get that almost surely, there exist infinitely many integers $n$ such that $\left|D_n^{\omega}\right|>\gamma\Phi(q^n)$. If $q>0$ is taken large enough ($\kappa,\gamma$ and $c_2$ are not depending on $q$), we can use the inequality (2.1) to prove that almost surely:
$$ \left|X_{q^n}^{\omega}\right|\geq \gamma \Phi(q^n)-2c_+q^{-1/2}\Phi(q^n)>\frac{\gamma}{2}\Phi(q^n)$$
for infinitely many $n$, which is the claim.
\begin{flushright}
$\ $
\end{flushright}

$ $\\
\textbf{Remark 1:} In order to bound the sum in (2.2) from below, Gaussian bounds were not sufficient. Without the volume growth property, the annulus could contain only few elements. Even if the exponential term is not too small (typically of order $n^{-s}$ for $s$ small), the term $t_n^{-d/2}$ (which corresponds to $t^{-d/2}$ for the Brownian motion) could be very small and make the series become summable. The cardinality of the annulus was critical in order to balance out this term.

$ $\\ 
\textbf{Remark 2:} our goal was to obtain a result in ${L}^1$ norm. Unfortunately, the
natural distance on graphs is the chemical distance $d_{\omega}(.,.)$. In the bound from below, this  does not create any trouble
because of the trivial inequality $\left|x\right|\leq d_{\omega}(0,x)$. But it could happen that the chemical distance is much bigger
than the $L^1$ norm. The proof of Theorem 2.18 in Barlow \cite{Barlow} precisely deals with this issue
thanks to a result by Antal and Pisztora \cite {AntalPisztora} that shows that the chemical
distance on $\mathcal{C}$ and the ${L}^1$ norm are not that different on a supercritical percolation cluster.

\section{Zero-one Law for the blind random walk} 

In the present section, we will consider discrete time random walks. We first introduce the two random walks we
will use. Then we recall an ergodicity result proved in \cite{BergerBiskup} and we derive the Zero-one Law. 
Our proofs are rather direct applications of the ergodicity statement of \cite{BergerBiskup}.

For each $x\in \mathbb{Z}^d$, let $\tau_x$ be the shift from $\Omega$ in $\Omega$ defined by: $(\tau_x\omega)_y=\omega_{y+x}$. For each $\omega$, let $Y_n^{\omega}$ be the simple random walk (called \textbf{blind random walk}) on $\mathcal{C}$ started at the origin. At each unit of time, the walk picks a neighbor at random and if the corresponding edge is occupied, the walk moves to this neighbor. Otherwise, it does not move. This random walk may seem less natural than the random walk that chooses randomly one of the accessible neighbors and jumps to it, but this blind random walk preserves the uniform measure on $\mathcal {C}$, so that the stationary measure on the environment as seen from the walker turns out to be simpler.  

It is well known (cf Kipnis and Varadhan \cite{KipnisVaradhan}) that the Markov chain $(Y_n^{\omega})_{n\geq 0}$ induces a Markov chain on $\Omega$ (the so-called \textbf{Markov chain on the environment}), that can be interpreted as the trajectory of "environment viewed from the perspective of the walk". It is defined as 
$$\omega_n ( \cdot) = \omega ( \cdot + Y_n^{\omega})=\tau_{Y_n^{\omega}}\omega(.).$$
One can describe the chain $(\omega_n)$ as follows. At each step $n$, one chooses one of the $2d$ neighbors of the origin at
random and calls it $e$. If the corresponding edge is closed for $\omega_n$, then $\omega_{n+1} = \omega_n$, otherwise $\omega_{n+1} (\cdot) = \tau_e \circ \omega_n$, where $\tau_e \circ \omega (\cdot) = \omega ( \cdot - e)$.    
 
It is straightforward to check that the probability measure $\tilde{\mathbb{P}}_p$ is a reversible and therefore stationary measure for the Markov chain $(\omega_n)$.
This allows us to extend our probability space to $\Xi= \Omega^{\mathbb{Z}}$ (endowed with the product $\sigma$-algebra $\mathcal {H}= \mathcal{F}^{\otimes\mathbb{Z}}$) and to define $\omega_n$ also for negative $n$'s in such a way that
that the family $(\omega_n, n \in \mathbb {Z})$ is stationary. Let $\mu$ denote the probability measure associated to the
Markov chain.
 
Note that under the measure $\mu$, and for all $n \in \mathbb {Z}$, the law of $(\omega_n, \omega_{n+1}, \ldots)$ is identical
to that of $(\omega_0, \omega_1, \ldots) $. On the other hand, the marginal law of $\omega_0$ (still under $\mu$) is $\tilde{\mathbb{P}}_p$. One then defines $T:\Xi \rightarrow \Xi,\bar{\omega}\mapsto T\bar{\omega}$ to be $(T\bar{\omega})_n=\bar{\omega}_{n+1}$. 

\begin{theorem}(Berger, Biskup, \cite {BergerBiskup})
\label{ergodicity} 
$T$ is ergodic with respect to $\mu$. In other words, for all $A\in \mathcal{H}$, if $T^{-1}(A) =A$, then $\mu(A)$ is equal to 0 or 1.
\end{theorem}

We refer to the paper of Berger and Biskup \cite{BergerBiskup} for proofs. Define for every $a>0$ and $\omega$ the event: $$A_{\omega}(a)=\left\{\limsup_{n\rightarrow \infty}\frac{\left|Y_n^{\omega}\right|}{\Phi(n)}>a\right\}.$$
Let now state and prove a consequence of this ergodicity for our law of the iterated logarithm: 

\begin{corollary}
\label{loi du tout ou rien} \textbf{(Zero-one Law)} Let $a\geq 0$. The probability that
$$B_a=\left\{\mathbb{P}_{\omega}^x\text{ a.s. }A_{\omega}(a)\text{ holds for all }x\in \mathcal{C}\right\}$$ is equal to $0$ or to $1$.
\end{corollary}

\textbf{Proof:} 
Our goal is to use the ergodicity of the environment and to note that the considered event corresponds to a $T$-invariant set in $\Xi$. 
Let $a \geq 0$ and define the function $F$ on $\Omega$ by:
$$F(\omega)=\mathbb{P}_{\omega}^0(A_{\omega}(a))$$
This function is well-defined and measurable. Let fix the environment $\omega$ for a little while and
denote $\omega_n=\tau_{Y_n^{\omega}}\omega$. We claim that $(F(\omega_n))_n$ is a martingale with respect to the
filtration $\mathcal{F}_n$ associated to the process $Y_n^{\omega}$. Indeed, the Markov property yields
$$F(\omega_n)=\mathbb{P}_{\omega}^0(A_{\omega}(a)\left|\mathcal{F}_n\right.).$$
This martingale is bounded and therefore converges almost surely as $n \to \infty$.
 Moreover, it converges to the indicator function 
of $A_{\omega}(a)$ because this event is clearly in
$\mathcal{F}_{\infty}=\sigma (\cup_{n\geq 0}\mathcal{F}_n)$. 

By taking the Cesaro mean and then integrating it with respect to $\omega$ (and using the fact that the probabilities are bounded by $1$), we get that
$$
 \tilde{\mathbb{E}}_0\left[
 \mathbb{P}_{\omega}^0 \left( \left| \lim_{N\rightarrow \infty} \frac{1}{N} \sum_{n=0}^{N-1}F(\omega_n) - 1_{A_{\omega}(a)} \right| \right) \right] = 0 $$

On the other hand, $F$ can be viewed as a measurable function on $\Xi$. The ergodicity of $\mu$ implies that for $\mu$ almost every
$\bar{\omega}$:
$$\lim_{N\rightarrow \infty} \frac{1}{N} \sum_{n=1}^{N}F(\bar{\omega}_n)= \int Fd\mu$$
Let recall that $\bar{\omega}$ has same law under $\mu$ as $(\omega_n)_n$ under
$\tilde{\mathbb{E}}_0[\mathbb{P}_{\omega}^0(.)]$. We deduce that the limit 
$1_{A_{\omega}(a)}$ is (up to a set of zero measure) constant. Since it is an indicator function, this means that 
either the corresponding event is almost surely true, or almost surely wrong.

\section {The Law of the Iterated Logarithm}

We can now derive the Law of the Iterated Logarithm. Let first note that the previous corollary immediately implies that for a fixed $p > p_c$, there exists a constant $c'(p,d) \in [0, \infty]$ such that for almost
every environment $\omega$, the blind random walk satisfies 
$$\limsup_{n\rightarrow \infty} \frac{\left|Y_n^{\omega}\right|_1}{\Phi(n)}=c'(p,d)$$
almost surely (just choose $c'(p,d)$ to be the supremum of the set of $a$'s such that the event $B_a$ is almost surely satisfied).

Our next goal is to show that the time scales for the two random walks are comparable.  Let $\omega \in \Omega_0$, define the real random variable $(T_n^{\omega})_n$ by $T_0^{\omega}=0$ and
$$T_{n+1}^{\omega}=\inf\left\{t>T_n^{\omega}, X_t^{\omega} \neq X_{T_n^{\omega}}^{\omega}\right\}.$$
Clearly, the Law of Large Numbers implies that for all $\omega \in \Omega_{0}$,
$T_{n+1}^{\omega} \sim n $ almost surely.

Let $\omega \in \Omega_0$, define in the same way the random variable $(U_n^{\omega})$ by $U_0^{\omega}=0$ and 
$$U_{n+1}^{\omega}=\inf\left\{p>U_n^{\omega}, Y_p^{\omega} \neq Y_{U_n^{\omega}}^{\omega}\right\}.$$
The $(U_{n+1}^{\omega}-U_n^{\omega})_n$ are not i.i.d. anymore. Conditionally on the environment and on the past up to the
$n$-th jump of $Y^{\omega}$, the law of $U_{n+1}^{\omega}- U_n^{\omega}$ is geometric and depends on the number $I(n)$ of incoming
open edges at $Y^{\omega}{U_n^{\omega}}$ (its mean is some function $f(I(n))$).

Ergodicity ensures that almost surely and for each $k \le 2d$, 
$$ \frac {1}{n} \sum_{j=1}^n 1_{ I(j) = k } \to i(k)$$
where $i(k)$ denotes the $\mu$-probability that $\omega_0$ has $k$ incoming open edges at the origin.  

Using the Law of Large Numbers for sums of independent geometric random variables of mean $f(k)$ for each $k$, we get
readily that for almost all $\omega\in \Omega_0$, 
$$\mathbb{P}_{\omega}^0\text{ a.s. }U_{n+1}^\omega / n  \to \sum_{k=1}^{2d} i(k) f(k)=\alpha_p^{-1}.$$
This last quantity is clearly positive and finite. 

We can now conclude the proof of the Law of the Iterated Logarithm for the continuous time random walk. 

\medbreak

\textbf{Proof of Theorem \ref{law of the iterated logarithm}:} Consider the natural coupling for which $X_t^{\omega}$ and
$Y_n^{\omega}$ have the same trajectories. More precisely, if we consider the \textbf{myopic random walk} $(Z_n^{\omega})_n$
that jumps at each time, choosing uniformly a neighbor, defined on a probability space $(\Omega_{\omega},\mathcal{F}_{\omega},\mathbb{P}^0_{\omega})$. Assume there exists an independent family $(T_i)_{i\in\mathbb{Z}_+}$ of iid exponential mean time 1 random variables and $(S^{\omega}_x)_{x\in \mathbb{Z}^d}$ an independent family of independent random variables such that $S_x^{\omega}$ is a geometrical of parameter ${n_x^{\omega}}/{(2d)}$ where $n_x^{\omega}$ is the number of adjacent open edges of $x$ for the configuration $\omega$. 

Define $T_p^{\omega}=\sum_{k=0}^{p-1}T_i$ and $n^{\omega}(t)=\sup\left\{p, T_p^{\omega} \leq t\right\}$. Then we can write the continuous time random walk as follows
$$X_t^{\omega}=Z_{n^{\omega}(t)}^{\omega}\ \ \forall t \geq 0.$$
Now, consider $U_p^{\omega}=\sum_{k=0}^{p-1}S_{Z_k^{\omega}}^{\omega}$ and $m^{\omega}(n)=\sup\left\{p, U_p^{\omega} \leq p\right\}$. Then we can write the blind random walk as follow:
$$Y_n^{\omega}=Z_{m^{\omega}(n)}^{\omega}\ \ \forall n \geq 0.$$
Because of the estimates of the time-scales of our two walks, we get that
$$
\limsup_{t \rightarrow \infty}\frac{\left|X_t^{\omega}\right|}{\Phi(t)} =  \limsup_{t \rightarrow \infty}\frac{\left|Z_{n^{\omega}(t)}^{\omega}\right|}{\Phi(n^{\omega}(t))}=\limsup_{n \rightarrow \infty}\frac{\left|Z_{n}^{\omega}\right|}{\Phi(n)}$$
and that
$$
\limsup_{n \rightarrow \infty}\frac{\left|Y_n^{\omega}\right|}{\Phi(n)} =\limsup_{n \rightarrow \infty}\frac{\left|Z_{
m^{\omega}(n)}^{\omega}\right|}{\Phi(\alpha_p m^{\omega}(n))}=\frac{1}{\sqrt{\alpha_p}}\limsup_{n \rightarrow \infty}\frac{\left|Z_{n}^{\omega}\right|}{\Phi(n)}.$$
From these two equalities, we deduce that
$$\limsup_{t \rightarrow \infty}\frac{\left|X_t^{\omega}\right|}{\Phi(t)}=\frac{1}{\sqrt{\alpha_p}}\limsup_{n \rightarrow \infty}\frac{\left|Y_n^{\omega}\right|}{\Phi(n)}\text{ a.s.}.$$
The theorem follows readily.
\begin{flushright}
$\square$
\end{flushright}

Note that this also show that the Law of the Iterated Logarithm holds for the blind and the myopic random walks.

\medbreak

\textbf{Acknowledgements:} This paper was written during my stay at the University of British Columbia, I would like to thank M.T. Barlow,
who first taught me about this question, for his availability and the advice he gave me during my whole stay. I would also
like to thank W. Werner for his careful reading of this paper and his numerous suggestions.

\begin{flushright}

\footnotesize \textsc{Department of Mathematics}

\textsc{University of British Columbia}

\textsc{Vancouver, British Columbia, Canada}

\medbreak

\textsc{DMA, Ecole Normale Sup\'erieure}

\textsc{45 rue d'Ulm, 75230 Paris cedex 05, France}

\textsc{E-mail:} hugo.duminil@ens.fr\end{flushright}
\end{document}